\documentclass[12pt]{amsart}
\usepackage[T1]{fontenc}
\usepackage[spanish]{babel} 
\usepackage{amsmath}
\usepackage{amscd}
\usepackage{graphics}
\usepackage{latexsym}
\usepackage[all]{xy}
\usepackage{epigraph}

\makeatletter
\renewcommand\@makefntext[1]{\leftskip=2em\hskip-2em\@makefnmark#1}
\makeatother

\newcommand{\RR}{\mathbb{R}}
\newcommand{\QQ}{\mathbb{Q}}

\newcommand{\ZZ}{\mathbb{Z}}
\newcommand{\CC}{\mathbb{C}}
\newcommand{\PP}{\mathbb{P}}

\begin{document}
\title[Curves]{Curvas algebraicas y la pregunta de Halphen} 

\bigskip
\author[Lozano]{C\'esar Lozano Huerta}
 \address{Instituto de Matem\'aticas\\
 UNAM. Oaxaca de Ju\'arez, Oaxaca. M\'exico}
 
 \address{Department of Mathematics \\
 Harvard University. Cambridge, MA. USA}

\email{lozano@math.harvard.edu} 

\maketitle

\section*{Introducci\'on}~\label{intro}

\medskip\noindent
Jakob Steiner fue un matem\'atico suizo del siglo {\footnotesize{XIX}}, que en vida se propuso compilar y modernizar la geometr\'ia sint\'etica, conocida desde el tiempo de los griegos del siglo {\footnotesize{IV}} a.$\!$ de C. A su muerte, su testamento dotaba de ocho mil t\'aleros\footnote{Ocho a\~nos despu\'es de su muerte el t\'alero se remplaz\'o por el marco alem\'an a una tasa de cambio de $1$ t\'alero = $3$ marcos alemanes. Ocho mil t\'aleros actualmente ser\'ian seis mil d\'olares norteamericanos apr\'oximadamente.} de forma bianual al autor del mejor trabajo en geometr\'ia abordado sint\'eticamente; es decir, sin el uso de coordenadas. De este modo naci\'o el Premio Steiner que otorgaba la Universidad de Berl\'in. En 1882 este premio se dividi\'o entre Max Noether y Henri Halphen por su investigaci\'on sobre curvas algebraicas. Halphen, en un tratado de doscientas p\'aginas \cite{HAL}, aborda el problema sobre la existencia de curvas en el espacio proyectivo de dimensi\'on $3$, mismo que tiene continuidad hasta nuestros d\'ias y que nos servir\'a de gu\'ia para presentar las ideas del presente art\'iculo. Enunciaremos la pregunta central de \cite{HAL}, que llamaremos pregunta de Halphen, luego explicaremos la terminolog\'ia y la geometr\'ia del problema. La pregunta, escrita en un lenguaje moderno, dice:

\medskip
\begin{center}
\textit{?`Qu\'e pares de n\'umeros $(d,g)$ ocurren como el grado y g\'enero de una curva algebraica suave en el espacio proyectivo de dimensi\'on $3$? }
\end{center}

\medskip\noindent
Para apreciar la geometr\'ia detr\'as de esta pregunta, necesitamos
precisar, qu\'e es una curva algebraica y qu\'e son el grado y el g\'enero. Primero, contestaremos estas tres cuestiones con cierta generalidad, y as\'i responderemos burdamente la pregunta de Halphen. Al momento de la publicaci\'on de \cite{HAL}, no exist\'ia el concepto de curva algebraica abstracta. Los ge\'ometras estudiaban las curvas algebraicas encajadas en un espacio proyectivo, y buscaban clasificarlas; de ah\'i la relevancia de la pregunta central de \cite{HAL}. Al intentar clasificar las curvas encajadas, el g\'enero es un invariante intr\'inseco, es decir, no depende del encaje y por ende es importante saber calcularlo en ejemplos concretos. En este escrito calcularemos expl\'icitamente el g\'enero de curvas en el plano proyectivo y en ejemplos particulares de curvas en el espacio proyectivo.

\medskip\noindent
Pese a que la pregunta de Halphen se plante\'o, al menos, hace ciento treinta a\~nos, la abordaremos con el lenguaje y punto de vista de la teor\'ia de esquemas. No definiremos formalmente lo que es un esquema, pero uno de los objetivos de este texto es ilustrar un aspecto central de dicha teor\'ia en el caso de curvas, el cual es: estudiar la geometr\'ia de una curva o de una familia de ellas, examinando objetos puramente algebraicos. Por tanto, asumiremos que el lector est\'a familiarizado con los conceptos de anillo de polinomios y sus ideales.

\section{Para fijar ideas}

\medskip\noindent
El \'algebra lineal se ocupa de estudiar las caracter\'isticas del conjunto soluci\'on de una colecci\'on de ecuaciones lineales,
\begin{equation*}
\begin{aligned}
a_{11}x_0+\cdots +a_{1n}x_n&=0\\
&\vdots\\
a_{m1}x_0+\cdots +a_{mn}x_{n}&=0.
\end{aligned}
\end{equation*}
Aqu\'i, los coeficientes $a_{ij}$ pertenecen a un campo $\mathbb{K}$, y el conjunto soluci\'on es un espacio vectorial de dimensi\'on finita sobre $\mathbb{K}$.

\medskip\noindent
La geometr\'ia algebraica se ocupa del caso general: las caracter\'isticas del conjunto soluci\'on de una colecci\'on de ecuaciones como
\begin{equation}\label{DOS}
\begin{aligned}
f_1(x_0,\ldots, x_n)&=0 \\
&\vdots \\
f_m(x_0,\ldots, x_n)&=0,
\end{aligned}
\end{equation}
donde cada $f_r$ es un polinomio de $n$ variables de grado arbitrario con coeficientes en el campo $\mathbb{K}$. A dicho conjunto soluci\'on lo llamaremos conjunto algebraico y lo denotaremos por $V(f_1,\ldots , f_m)$.

\medskip\noindent
A diferencia del caso donde los polinomios tienen grado $1$, el conjunto soluci\'on de las ecuaciones en (\ref{DOS}) no es un espacio vectorial, y adem\'as es sensible a las caracter\'isticas algebraicas del campo $\mathbb{K}$ al que pertenecen los coeficientes. Por simplicidad, trabajaremos con el campo de los n\'umeros complejos $\CC$. No obstante, dado que $\CC$ es un espacio vectorial sobre $\RR$ de dimensi\'on $2$, usar n\'umeros complejos implica considerar espacios de dimensi\'on real $2$, $4$ y mayores. 

\medskip\noindent
Al estudiar las soluciones de ecuaciones como las que aparecen en (\ref{DOS}), las preguntas fundamentales que surgen son: ?`qu\'e tipo de conjunto forman las soluciones?, ?`es un conjunto finito?, si \'este no es finito ?`es un conjunto continuo?, ?`de  qu\'e dimensi\'on?
Nos enfocaremos en entender el conjunto soluci\'on de una colecci\'on de ecuaciones del tipo (\ref{DOS}) cuando \'este tenga dimensi\'on $1$ sobre los complejos; conjunto al que podr\'iamos llamar curva.
!`Cuidado!, no hemos dicho a\'un qu\'e es la dimensi\'on de un conjunto soluci\'on; lo diremos m\'as adelante usando la funci\'on de Hilbert. De momento, algo m\'as fundamental es comentar d\'onde habita el conjunto soluci\'on $V(f_1,\ldots , f_m)$. Dado que nos interesa aprovechar las herramientas de la geometr\'ia proyectiva, estudiaremos el caso cuando dicho conjunto es un subconjunto del espacio proyectivo de dimensi\'on $n$, el cual denotamos por $\PP^n$. Abundaremos sobre este espacio m\'as adelante y veremos que trabajar con \'el implica que s\'olo estudiaremos conjuntos algebraicos $V(f_1,\ldots , f_m)$, donde los polinomios $f_r$ son homog\'eneos\footnote{Por definici\'on, $f$ es homog\'eneo de grado $d$ si $f(\lambda x_0,\ldots , \lambda x_n)=\lambda^d f(x_0,\ldots ,x_n)$ para todo $\lambda \in \CC\backslash \{0\}$.}.

\medskip\noindent
Resumiendo, estudiaremos el conjunto $$V(f_1,\ldots , f_m)\subset \PP^n,$$ cuando \'este tiene dimensi\'on $1$, y cada $f_r\in \CC[x_0,\ldots , x_n]$ es un polinomio homog\'eneo. La ventaja de estudiar curvas contenidas en $\PP^n$ es que las podemos caracterizar de manera burda, usando s\'olo dos n\'umeros enteros: el grado y el g\'enero.

\medskip\noindent

\section{Curvas en el plano proyectivo $\mathbb{P}^2$}\label{UNO}

\medskip\noindent
Como primer ejemplo concreto estudiaremos curvas planas. Textos introductorios sobre el tema son \cite{KIRWAN,REID}. En esta secci\'on, calcularemos expl\'icitamente su g\'enero y una versi\'on de espacio tangente. Estas curvas se llaman planas porque est\'an contenidas en el plano proyectivo complejo $\PP^2$, el cual tiene dimensi\'on $2$ sobre los complejos, pero !`dimensi\'on $4$ sobre los reales! Este espacio posee coordenadas globales $[x:y:z]$ que est\'an sujetas a la condici\'on $[x:y:z]=[\lambda x:\lambda y:\lambda z]$ para todo n\'umero complejo $\lambda\ne 0$ y $(x,y,z)\ne(0,0,0)$. Por ejemplo, $[1:2:3]$ y $[2:4:6]$ son coordenadas del mismo punto en $\PP^2$. Espec\'ificamente, denotando $(x,y,z)\in \CC^3\backslash \{(0,0,0)\}$, se define $$\PP^2:=\{[x:y:z]\  |\  (x,y,z)=(\lambda x,\lambda y,\lambda z) \mbox{ y }\lambda\in \CC^*\}.$$ 

\medskip\noindent
Una curva plana est\'a definida por una ecuaci\'on en las coordenadas de $\PP^2$. Es decir, sea $f\in \CC[x,y,z]$ un polinomio de tres variables con coeficientes complejos, homog\'eneo de grado $d$. Definimos una \textit{curva algebraica plana} $C$ como 
$$C=\{[x:y:z]\ \in \PP^2\ |\ f(x,y,z)=0\}.$$
El grado de la curva $C$ se define como el grado del polinomio $f$. Escribimos $C=V(f)$ cuando deseamos hacer referencia expl\'icita al polinomio que define la curva. Observar que la propiedad $[x:y:z]=[\lambda x:\lambda y:\lambda z]$ para todo n\'umero complejo $\lambda\ne 0$, hace que solo tenga sentido hablar de un conjunto algebraico $V(f)$ cuando $f$ es homog\'eneo.

\medskip \noindent 
El concepto de suavidad es f\'acil de imaginar, y para curvas planas sencillo de definir. Una curva plana $C=V(f)$ es suave en un punto $p\in C$ si la \textit{recta tangente} a $C$ en $p=[p_1:p_2:p_3]$, denotada como $T_pC$, est\'a bien definida. \'Esta \'ultima se define como: $$\frac{df}{dx}(p)(x-p_1)+\frac{df}{dy}(p)(y-p_2)+\frac{df}{dz}(p)(z-p_3)=0.$$

\begin{center}
\begin{figure}[htb]
\resizebox{1\textwidth}{!}{\includegraphics{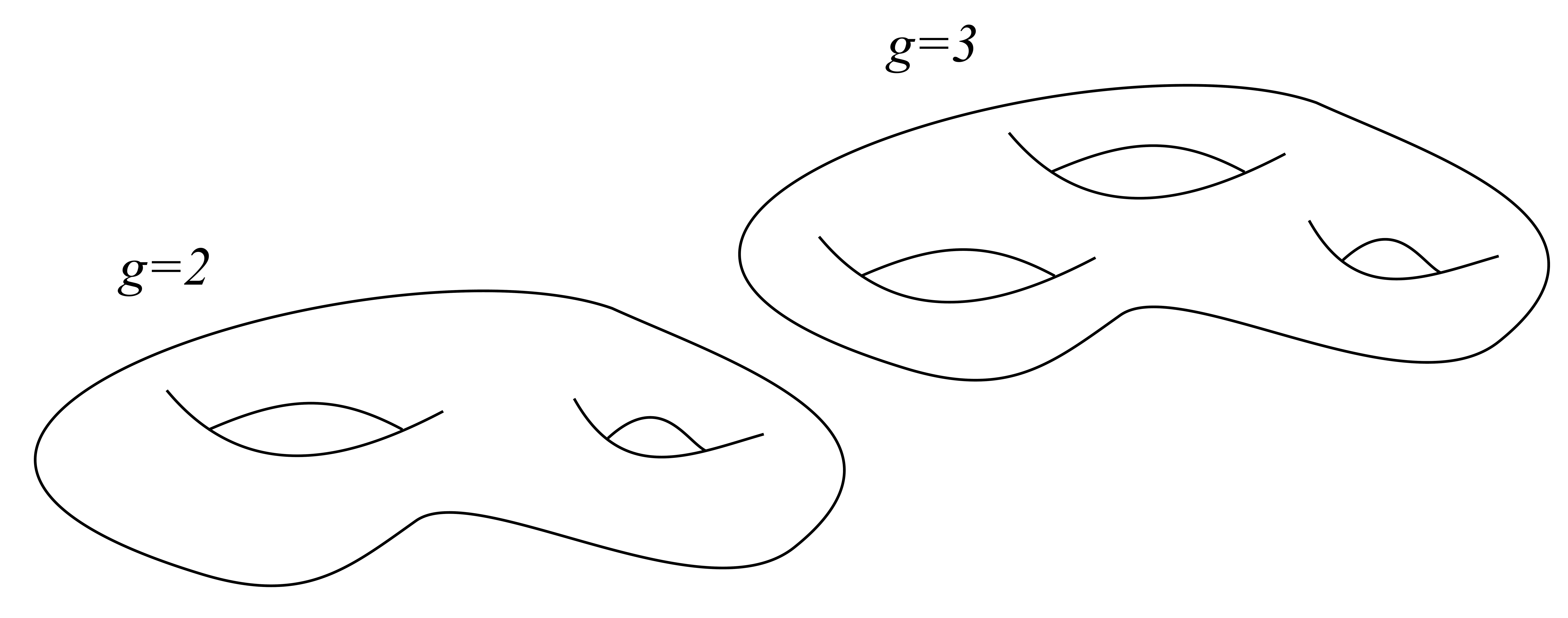}}
\caption{Curvas suaves de g\'enero $2$ y $3$.}\label{GENUS}
\end{figure}
\end{center}

\medskip\noindent 
Antes de definir formalmente el g\'enero, haremos un comentario sobre lo que representa geom\'etricamente. Las curvas planas suaves, y de hecho todas las curvas suaves en $\PP^n$, son compactas, conexas, de dimensi\'on $2$ sobre $\RR$.\footnote{Por esto en geometr\'ia compleja a las curvas algebraicas suaves se les llama superficies de Riemann.} Por lo tanto, topol\'ogicamente, s\'olo tienen un invariante: el g\'enero.
El g\'enero de una curva $C$ es un n\'umero dif\'icil de definir, pero muy f\'acil de imaginar. La figura \ref{GENUS} muestra dos curvas de g\'enero $2$ y $3$.

\noindent 
La figura \ref{GENUS} deja clara la idea de lo que significa (topol\'ogicamente) una curva suave de g\'enero $2$ y $3$ \textemdash \ \!salvavidas para $2$ y $3$ personas, respectivamente\textemdash. De esta idea intuitiva, podemos concluir que el g\'enero de una curva suave tiene que ser al menos $0$; como la esfera, la cual es, en efecto, la curva de menor g\'enero.

\medskip \noindent 
Lo que sigue es la parte t\'ecnica, pues definiremos formalmente el g\'enero de una curva plana. El lector notar\'a que este n\'umero, el cual posee informaci\'on geom\'etrica, se definir\'a usando propiedades de los polinomios. Sea $f\in R=\CC[x,y,z]$ un polinomio de tres variables homog\'eneo de grado $d$. Consideremos el anillo cociente $R/\langle f \rangle$, donde $\langle f \rangle$ denota el ideal generado por $f$. Observemos que los elementos del anillo $R$, al igual que los del anillo cociente $R/\langle f \rangle$, los podemos distinguir por su grado. La funci\'on que contabiliza los elementos que hay en cada grado tiene nombre propio: funci\'on de Hilbert.  Por ejemplo, en el anillo $R$ hay tres monomios de grado $1$, $R_1=(x,y,z)$, y dado que todos los polinomios de grado $1$ son combinaciones lineales en $x,y$ y $z$, entonces la funci\'on de Hilbert  toma el valor $H_R(1)=\mathrm{dim}_{\CC}\ \! R_1=3$. Similarmente en grado $2$, el anillo $R$ contiene seis monomios $R_2=(x^2,y^2,z^2,xy,xz,yz)$ y todos los dem\'as son combinaciones lineales de \'estos, luego $H_R(2)= 6$. En general, $H_R(m)$ es igual al n\'umero de monomios en grado $m$ en $R$, que en este caso es $H_R(m)=\mathrm{dim}_{\CC}\ \! R_m=\frac{1}{2}(m+2)(m+1)$.
Si volvemos a contar, ahora para el cociente $R/\langle  f \rangle$, la funci\'on que nos dice la dimensi\'on de cada componente homog\'enea de grado $m$ del anillo cociente es:
\begin{equation}\label{HF}
H_C(m):=\mathrm{dim }_{\CC}(R/\langle f \rangle)_m. 
\end{equation}
A la funci\'on $H_C(m)$ le llamaremos la funci\'on de Hilbert de la curva $C=V(f)\subset \PP^2$. 

\medskip\noindent
Conocer el comportamiento de $H_C(m)$ es un objetivo central de este art\'iculo y un resultado debido a David Hilbert lo estima: existe un polinomio $P_C(m)$, que toma los valores de $H_C(m)$ para valores grandes de $m$ \cite[p.~51]{HART}. Es decir, la funci\'on $H_C(m)$ se comporta asint\'oticamente como un polinomio; el cual se sigue que es \'unico. Dicho polinomio tiene informaci\'on geom\'etrica de la curva $C$; como lo vemos a continuaci\'on.

\medskip\noindent
Sea $P_C(m)$ el polinomio asociado a la funci\'on de Hilbert $H_C(m)$ de la curva $C=V(f)\subset \PP^2$, descrito en el p\'arrafo anterior. Definimos el \textit{g\'enero} de $C$ como $g(C):=1-P_C(0)$. Un poco m\'as adelante calcularemos $P_C(m)$ explicitamente. !`Ojo! pese a que $H_C(m)=P_C(m)$ para valores grandes de $m$ y por tanto son valores positivos, el valor de $P_C(0)$ podr\'ia ser negativo; de hecho es negativo en la mayor\'ia de los casos cuando $C$ es suave.

\medskip\noindent
No es obvio que la noci\'on de g\'enero que expresa la figura 1 y el n\'umero que acabamos de definir coinciden cuando la curva es suave. El lector puede (y debe) consultar \cite[p.~131]{MUM} para una prueba geom\'etrica de este hecho. Es tan poco obvio que el autor de esta \'ultima referencia, David Mumford, comenta que la igualdad de estos dos n\'umeros es una consecuencia del Teorema de Hirzebruch-Riemann-Roch. Con este comentario, Mumford enfatiza que la igualdad entre los g\'eneros yace en la intersecci\'on de tres ramas de las matem\'aticas: geometr\'ia, an\'alisis y topolog\'ia. 
Dado que en este art\'iculo deseamos enfatizar las ideas alrededor de la pregunta de Halphen no abordaremos la prueba de este hecho en esta ocasi\'on. 

\medskip\noindent
Una de las ventajas de la definici\'on $g(C)=1-P_C(0)$ es que no depende del campo $\CC$, ni de la suavidad de $C$. S\'olo depende de $H_C$, la funci\'on de Hilbert, la cual calculamos a continuaci\'on en el caso de curvas planas.

\medskip\noindent
Hemos calculado antes la funci\'on $H_R(m)$, con $R=\CC[x,y,z]$, y de manera similar podemos calcular la funci\'on $H_f(m)$. Esta \'ultima es la funci\'on que contabiliza la dimensi\'on del espacio de polinomios de grado $m$ en el ideal $\langle f \rangle$. Por definici\'on, el ideal $\langle f \rangle$ se compone de todos los m\'ultiplos de $f$ y por tanto, un elemento de grado $d+1$ en $\langle f \rangle$ es simplemente $lf$, donde $l$ es un polinomio lineal. Por lo tanto, la dimensi\'on de la componente de grado $m\ge d+1$ en $\langle f \rangle$ es igual al n\'umero de monomios de grado $m-d$. Por \'unica vez, usaremos la sucesi\'on de $R$-m\'odulos
\begin{equation*} \label{SEQ}
0\rightarrow \langle f \rangle \overset{f}{\rightarrow} R\rightarrow 
R/\langle f \rangle\rightarrow 0 ,\end{equation*}
pues \'esta implica que $H_R=H_{f}+H_{C}$ y por tanto, para $m\ge d+1$, obtenemos :
\begin{align*}
H_C(m) =& \{\mbox{monomios de grado $m$ en }R\}-\{\mbox{polinomios de grado $m$ en el ideal }\langle f \rangle\}\\
=&\binom{2+m}{2}-\binom{2+m-d}{2}.
\end{align*}

\medskip\noindent
Observando este c\'alculo inmediatamente nos damos cuenta de que $H_C(m)$ es un polinomio si $m\ge d+1$. Es decir, este c\'alculo nos ha dado $P_C(m)$. Concluimos que el g\'enero $g$ de una curva plana $C$ de grado $d$ es: 
\begin{equation}\label{GG} 
g=1-P_C(0)=\frac{(d-1)(d-2)}{2}.
\end{equation}

\medskip\noindent 
De (\ref{GG}) se sigue: !`no existen curvas planas (suaves o no) de g\'enero $2$! M\'as a\'un, si fijamos el grado de la curva, entonces su g\'enero est\'a determinado. No hay misterio sobre los pares $(d,g)$ que ocurren como el grado y g\'enero de las curvas planas. Esto contesta satisfactoriamente la pregunta de Halphen en el plano.
\bigskip

\section{Curvas en el espacio proyectivo $\PP^3$}~\label{SECTION2}

\noindent
En esta secci\'on abordaremos la pregunta de Halphen para curvas contenidas en $\PP^3$. El espacio proyectivo $\PP^3$ tiene coordenadas globales $[x:y:z:w]$ y a diferencia de $\PP^2$, todas las curvas algebraicas suaves se pueden encajar en \'el \cite[p.~ 310]{HART}. En t\'erminos m\'as precisos, $\PP^3$ se define como
$$\PP^3:=\{[x:y:z:w]\  |\  (x,y,z,w)=(\lambda x,\lambda y,\lambda z,\lambda w) \mbox{ y }\lambda\in \CC^*\}.$$ 

\medskip\noindent
Este espacio de dimensi\'on $3$ sobre los complejos, !`tiene dimensi\'on $6$ sobre los reales! En $\PP^3$ una curva est\'a definida al menos por dos ecuaciones y por tanto no podemos definir el grado como lo hicimos para curvas planas. Adem\'as, dos curvas del mismo grado pueden tener g\'eneros distintos (adelante citaremos un ejemplo), y por ende en $\PP^3$ no existe una ecuaci\'on como (\ref{GG}). La generalizaci\'on de estos invariantes yace en la funcion de Hilbert.

\medskip\noindent
En este momento nos enfrentamos a las preguntas: ?`qu\'e es una curva en $\PP^3$? y ?`c\'omo definimos su grado y g\'enero? Responderemos ambas de un solo brochazo con la ayuda de la funci\'on de Hilbert; lo que significa contestaremos examinando un objeto algebraico.

\medskip\noindent
Lo que sigue es t\'ecnico pues definiremos qu\'e es una curva algebraica en $\PP^3$. Nos auxiliaremos de la funci\'on de Hilbert generalizando la defini\'on (\ref{HF}). Dados varios polinomios homog\'eneos (ecuaciones) de grados arbitrarios en cuatro variables, $\{ f_1,\ldots ,f_k \}$, consideramos el ideal que generan $\langle f_1,\ldots ,f_k \rangle\subset R=\CC[x,y,z,w]$ y definimos la funci\'on de Hilbert de $V(f_1,\ldots, f_k)$ como: $$H(m):=\mathrm{dim}_{\CC} \left( R/\langle f_1,\ldots ,f_k \rangle\right)_m.$$

\medskip\noindent
Hemos dicho ya en la secci\'on anterior que conocer el comportamiento de $H_C(m)$ es un objetivo central de este art\'iculo, y que dicho comportamiento lo estima un resultado de David Hilbert \cite[p.~ 51]{HART}: existe un polinomio $P(m)\in \QQ[m]$, tal que $P(m)=H(m)$ para valores grandes de $m$. En la literatura, a este polinomio se le llama polinomio de Hilbert y sus coeficientes y grado tienen informaci\'on geom\'etrica importante del conjunto algebraico $V(f_1,\ldots, f_k)$ como lo indica el siguiente p\'arrafo. 

\medskip\noindent
Diremos que el conjunto algebraico $$V(f_1,\ldots, f_k)=\{\mathrm{\mbox{\boldmath$p$}} \in \PP^3|\ f_1(\mathrm{\mbox{\boldmath$p$}})=0,\ldots ,f_k(\mathrm{\mbox{\boldmath$p$}})=0\}$$ define una curva $C\subset \PP^3$, si su polinomio de Hilbert es de la forma $P_C(m)=Am+B$, con $A,B\in \ZZ$, y adem\'as es conexo.
El grado de $C$ se define como el n\'umero $A$ y su \textit{g\'enero} como $g(C):=1-B$. Por lo tanto, el polinomio de Hilbert de una curva $C$ de grado $d$ y g\'enero $g$ se escribe como 
\begin{equation}\label{PH} 
P_C(m)=dm-g+1.
\end{equation}

\medskip\noindent
El lector debe advertir que el grado de $P_C(m)$ nos dice qu\'e tan r\'apido crece la funci\'on $H_C(m)$ conforme $m$ crece y que estamos definiendo la dimens\'on de $C$ como este n\'umero. El grado del polinomio de Hilbert de un conjunto algebraico (es decir, c\'omo se comporta asint\'oticamente la funci\'on de Hilbert) define su dimensi\'on en un contexto m\'as general \cite[p.~ 51]{HART}.

\medskip\noindent
No definiremos suavidad de una curva en $\PP^3$, sin embargo este concepto formaliza la siguiente idea: una curva $C$ es suave en $p\in C$ si la recta tangente a $C$ en $p$, denotada como $T_pC$, est\'a bien definida y tiene dimensi\'on compleja $1$ \cite[p.~6]{MUM}, ver figura \ref{SUAVE}. La curva $C$ se dice suave si es suave para todo $p\in C$. 

\begin{center}
\begin{figure}[htb]\label{SUAVE}
\resizebox{.6\textwidth}{!}{\includegraphics{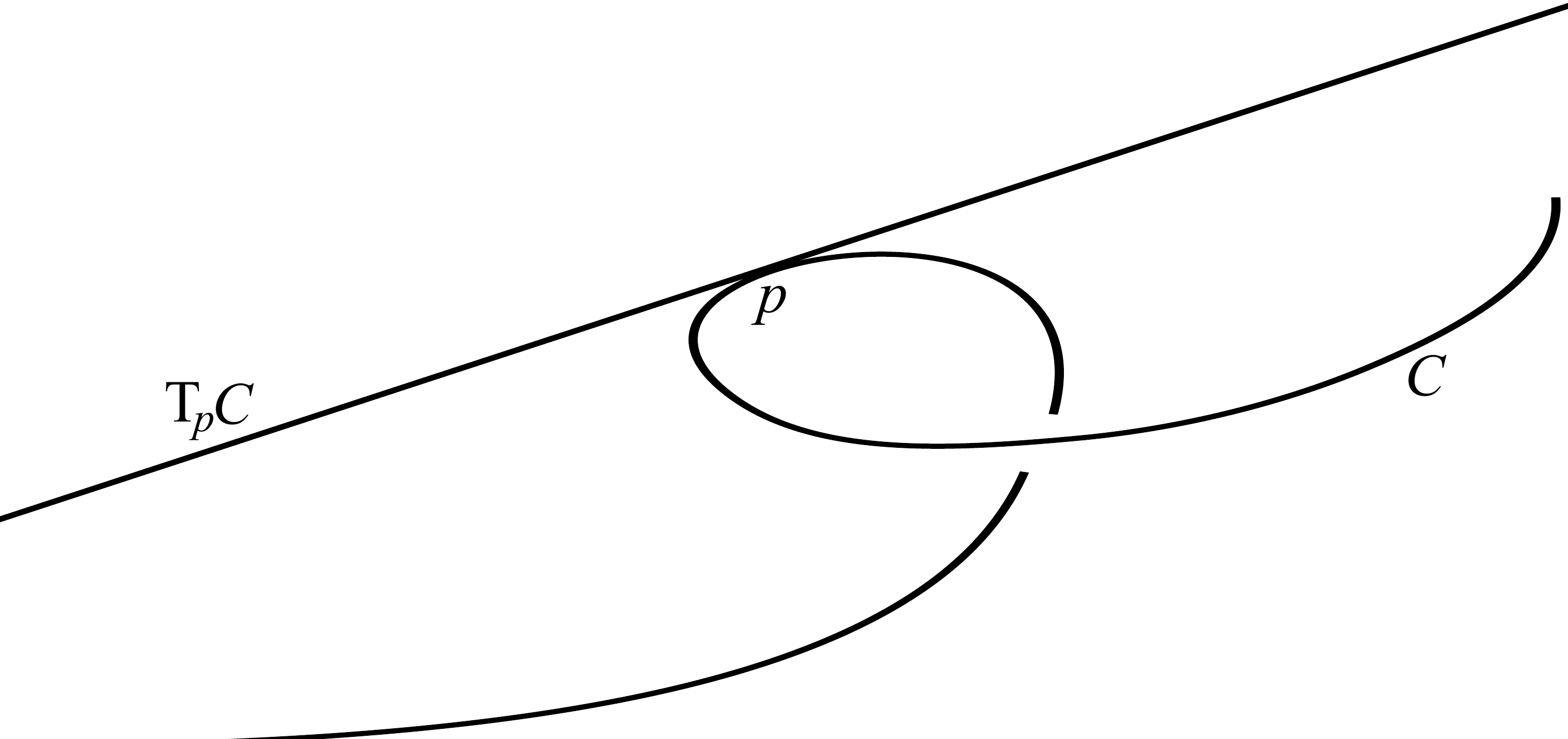}}
\caption{Recta tangente $T_pC$ a la curva $C$ en $p$.}\label{SUAVE}
\end{figure}
\end{center}

\medskip\noindent
El lector puede tambi\'en notar que en el caso $C\subset \PP^2$, la funci\'on de Hilbert $H_C(m)$ y polinomio de Hilbert $P_C(m)$ coinciden para todos los valores de $m\ge d+1$. En general, $H_C(m)\ne P_C(m)$ para valores peque\~nos de $m$. En el caso $C\subset \PP^3$, analizando la diferencia entre $H_C(m)$ y $P_C(m)$ para todos los valores de $m$ se obtendr\'a una cota para el g\'enero de $C$, escrita en (\ref{BOUND}), llamada \textit{cota de Castelnuovo}. Esto no resolver\'a la pregunta de Halphen, pero ser\'a un avance significativo hacia su soluci\'on como veremos en la siguiente secci\'on.

\medskip\noindent 
?`C\'omo se comporta el conjunto $V(f_1,\ldots ,f_k)$ si omitimos la condici\'on de suavidad? Veamos un ejemplo, el conjunto $C_0=V(y^2-zx,yw-z^2)$ tiene polinomio de Hilbert $P(m)=4m$, lo cual implica que es un conjunto algebraico de dimensi\'on $1$, cuyo g\'enero es $1$. Esta curva $C_0$ contiene el punto $p=[1:0:0:0]$ el cual no es suave: la recta tangente a $C_0$ en $p$ no esta bien definida pues existen dos curvas suaves, $L$ y $C$, contenidas en $C_0$ que pasan por $p$. M\'as a\'un, estas dos curvas forman el conjunto $C_0$. Por lo tanto, imponer la propiedad de suavidad a $V(f_1,\ldots ,f_k)$ evita que trabajemos con conjuntos como el anterior, que constan de varias piezas suaves. Abundaremos sobre el punto singular $p$ y las curvas $C_0, C$ y $L$, dos p\'arrafos adelante.

\medskip\noindent
Hemos definido el g\'enero independientemente de la suavidad de $V(f_1,\ldots ,f_k)$, y del campo $\CC$. Esta definici\'on adem\'as hace que el g\'enero (de hecho, el polinomio de Hilbert) sea constante en familias de curvas ``bien portadas'' \cite[p.~261]{HART}. Este \'ultimo adjetivo entrecomillado se formaliza en el concepto de \textit{familia plana}, del cual no hablaremos aqu\'i, pero que es importante en el estudio avanzado de curvas algebraicas \cite{MC}. 

\medskip\noindent
Veamos tres ejemplos. Sea $t$ un par\'ametro en el disco unitario $t\in \Delta=\{t\in \CC:\ |t|\le 1\}$.
\begin{equation}\label{RESIDUAL}
\begin{aligned}
    E&=V(w,zy^2-x^3+xz^2+z^3), & P_{E}(m)=3m+0 \\[3mm]
    C&=V(y^2-zx,yw-z^2,yz-xw), & P_{C}(m)=3m+1\\[3mm]
    C_t&=V(y^2-zx,yw-z^2+tx^2), & P_{C_t}(m)=4m+0.    
\end{aligned}
\end{equation}
Los polinomios de Hilbert listados aqu\'i nos dicen que $E$ tiene grado $3$ y g\'enero $1$, mientras $C$ tiene grado $3$ y g\'enero $0$. Por otro lado, para todo $t$, el conjunto $C_t$ tiene g\'enero $1$. Si $t=0$, es claro que $C\subset C_0$. M\'as a\'un, la curva $C$ junto con la l\'inea $L=V(y,z)$ forman $C_0$. El fen\'omeno que ocurre aqu\'i es que $C_t$ es suave de g\'enero $1$ si $t\ne 0$, y se ``quiebra'' en la uni\'on $C_0=C\cup L$, donde las componentes $C$ y $L$ son suaves y cada una tiene g\'enero $0$, a\'un cuando $C_0$ tiene g\'enero $1$. Los puntos singulares de $C_0$ son $C\cap L$; en particular, el punto $p=[1:0:0:0]$ es singular. El conjunto $C_t$ es un ejemplo de familia (plana) de curvas algebraicas en $\PP^3$.

\medskip\noindent
Resumiendo, el polinomio de Hilbert nos dice la dimensi\'on, g\'enero y grado del conjunto $V(f_1,\ldots ,f_k)$ sin importar si \'este es suave o no. Evitaremos trabajar con casos muy complicados o redundantes (como el ejemplo anterior), estudiando el caso suave. 

\medskip\noindent
A la luz de los ejemplos anteriores, reformulamos la pregunta original de Halphen: sea $C\subset \PP^3$ una curva suave de grado $d$, ?`cu\'ales son los n\'umeros $g$ que pueden ocurrir como el g\'enero de $C$? Las curvas $C$ y $E$ muestran que los pares $(d,g)=(3,0)$ y $(3,1)$ ocurren.

\medskip\noindent
Supongamos $C\subset \PP^3$ es una curva suave de grado $d$. Si $C$ yace en un plano, entonces la pregunta original de Halphen la contesta la ecuaci\'on (\ref{GG}). Esta f\'ormula, sin embargo, nos dice m\'as: el g\'enero $g(C)$ de una curva $C$ est\'a acotado, pues podemos proyectar dicha curva $C$ en un plano sin cambiar su grado \cite[p.~72]{MUM}, ver figura \ref{PROYECTAR}. Esto implica que
\begin{equation}\label{COTA} g(C)\le \frac{1}{2}(d-1)(d-2). \end{equation} Por lo tanto, para cualquier $d$ y $g(C)\ge \frac{1}{2}(d-1)(d-2)$, no existe curva en $\PP^3$ de este grado y g\'enero.

\begin{center}
\begin{figure}[htb]\label{PROYECTAR}
\resizebox{1\textwidth}{!}{\includegraphics{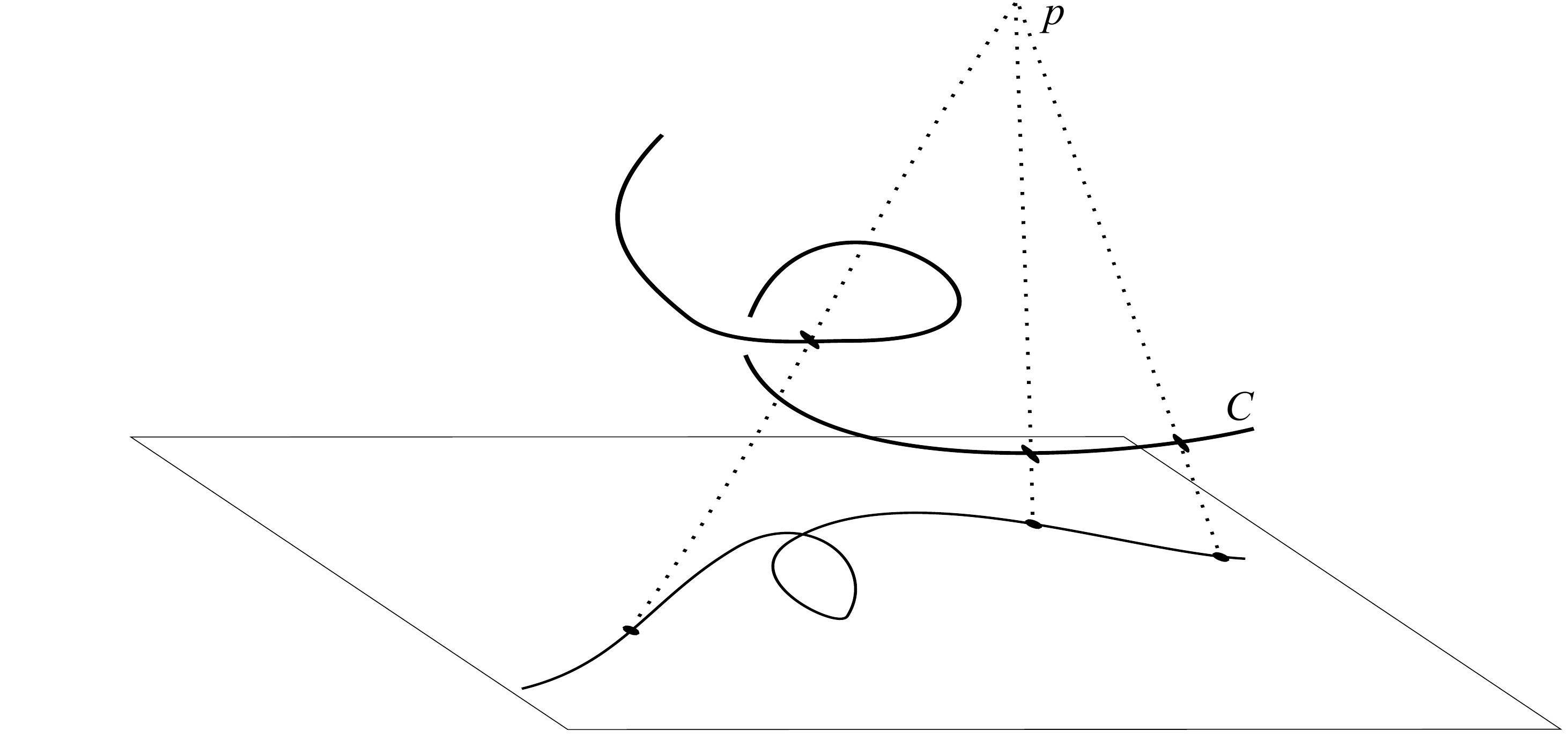}}
\caption{Proyecci\'on desde $p$ de una curva $C$ en un plano.}\label{PROYECTAR}
\end{figure}
\end{center}

\medskip\noindent
Halphen, y poco despu\'es Guido Castelnuovo\footnote{Guido Castelnuovo (1865-1952) naci\'o y creci\'o en Venecia. Se educ\'o y trabaj\'o, la mayor parte de su vida, en Italia. Tuvo estudiantes prominentes como O. Zariski y F. Enriques. Con el advenimiento del nazismo, Castelnuovo fue perseguido bajo leyes raciales fascistas y fue orillado a esconderse; al igual que muchos otros jud\'ios en Europa. Despu\'es de la guerra, Castelnuovo regres\'o a Italia y fue nombrado comisionado especial para reparar el da\~no a la ciencia italiana causado por el gobierno de veinte a\~nos de B. Mussolini. Muri\'o en Roma siendo senador de la Rep\'ublica Italiana.} \cite{CASTEL}, demostraron que si la curva $C$ no est\'a contenida en un plano, entonces $g(C)$ est\'a sujeto a una cota m\'as estricta. A saber,
\begin{equation}\label{BOUND}
  g(C) \le
  \begin{cases}
    \frac{1}{4}d^2-d+1 & d\text{ par }
    \\[3mm]
    \frac{1}{4}(d^2-1)-d+1 & d\text{ impar.}
  \end{cases}
\end{equation}

\medskip\noindent
Esto muestra que entre $0$ y la cota m\'as alta, dada por (\ref{COTA}), no todos los enteros ocurren. Es decir, entre $0$ y $\frac{1}{2}(d-1)(d-2)=\tfrac{1}{2}d-d+1$ existen ``huecos'' para los valores de $g$. Por ejemplo, se sigue de (\ref{BOUND}) que no existen curvas suaves de grado $4$ y g\'enero $2$, como tampoco existen curvas suaves de grado $5$ y g\'enero $3,4$ o $5$. 

\medskip\noindent
En 1889, Castelnuovo demostr\'o la desigualdad (\ref{BOUND}), ahora llamada cota de Castelnuovo, y su generalizaci\'on a curvas en $\PP^n$ con un argumento ingenioso que a continuaci\'on esbozaremos. Su estilo de hacer geometr\'ia algebraica influenciar\'ia a toda una generaci\'on de ge\'ometras italianos de principios del siglo {\footnotesize{XX}} \textemdash coloquialmente en estos d\'ias a dicha generaci\'on le llamamos ``la escuela italiana de geometr\'ia algebraica''\textemdash. Las ideas introducidas por \'el en \cite{CASTEL}, devinieron en lo que hoy se llama teor\'ia de Castelnuovo \cite{CP}.

\medskip\noindent
Esbocemos ahora el argumento de Castelnuovo para demostrar la desigualdad (\ref{BOUND}). Supongamos $C\subset \PP^3$ es una curva suave de grado $d$ y g\'enero $g$ que no est\'a contenida en ning\'un plano y est\'a definida por los polinomios $\{f_1,\ldots ,f_k\}$. Castelnuovo estima de dos maneras distintas el siguiente n\'umero$$H_C(m)=\mathrm{dim}_{\CC} \left( R/\langle f_1,\ldots ,f_k \rangle\right)_m.$$

\medskip\noindent
Una de esas estimaciones la conocemos para valores grandes de $m$: dicho n\'umero est\'a dado por el polinomio de Hilbert, $H_C(m)=P_C(m)=md-g+1$.\footnote{El teorema de Riemann-Roch, $h^0(m)-h^1(m)=md-g+1$, puede usarse tambi\'en en esta parte del argumento. En este art\'iculo, es m\'as natural usar el polinomio de Hilbert teniendo la ventaja de no depender del los n\'umeros complejos, ni de la suavidad de $C$.} Castelnuovo obtiene el teorema al estimar este mismo n\'umero --\ \!\!usando geometr\'ia proyectiva\ \!\!-- para todos los valores de $m$ y comparar con $P(m)$. Estas estimaciones, y la comparaci\'on entre ellas, se ven as\'i: $$md-g+1\ge r(r+2)+(m-r)d+1,$$ donde $d=2r+1$. De esta desigualdad obtenemos la cota deseada (\ref{BOUND}) para el g\'enero $g$. La parte derecha de esta desigual, y la desigualdad misma son los elementos de la prueba de Castelnuovo que estamos omitiendo. El lector interesado en leer la prueba completa puede consultar cualquiera de las siguientes referencias, \cite[p.~351]{HART} y \cite[p.~251]{GH}, cada una con enfoques ligeramente distintos.

\bigskip\noindent
Una curva con g\'enero m\'aximo, seg\'un la desigualdad (\ref{BOUND}), tiene nombre (o m\'as bien apellido) en la literatura: \textit{curva de Castelnuovo}. Dicha curva tiene la propiedad notable de estar contenida en una \'unica superficie cu\'adrica. Este hecho (no expl\'icito en el p\'arrafo anterior) brinda herramientas para estudiarlas, pues las curvas contenidas en cu\'adricas son bien conocidas. Por ejemplo, esta propiedad inmediatamente nos dice que una curva de Castelnuovo $C$ de grado $d=2k-1$, aparece en la intersecci\'on de una cu\'adrica $Q$ con una superficie $S$ de grado $k$, y que $Q\cap S=C\cup L$, donde $L$ es una l\'inea; la curva $C$ en (\ref{RESIDUAL}) es una curva de Castelnuovo para la cual $k=1$. Esta propiedad, y las ideas detr\'as de ella, ser\'a crucial para responder a la pregunta central de este escrito, como lo detallamos en la siguiente secci\'on.

\begin{center}
\begin{figure}[htb]\label{GRAFICA}
\resizebox{.8\textwidth}{!}{\includegraphics{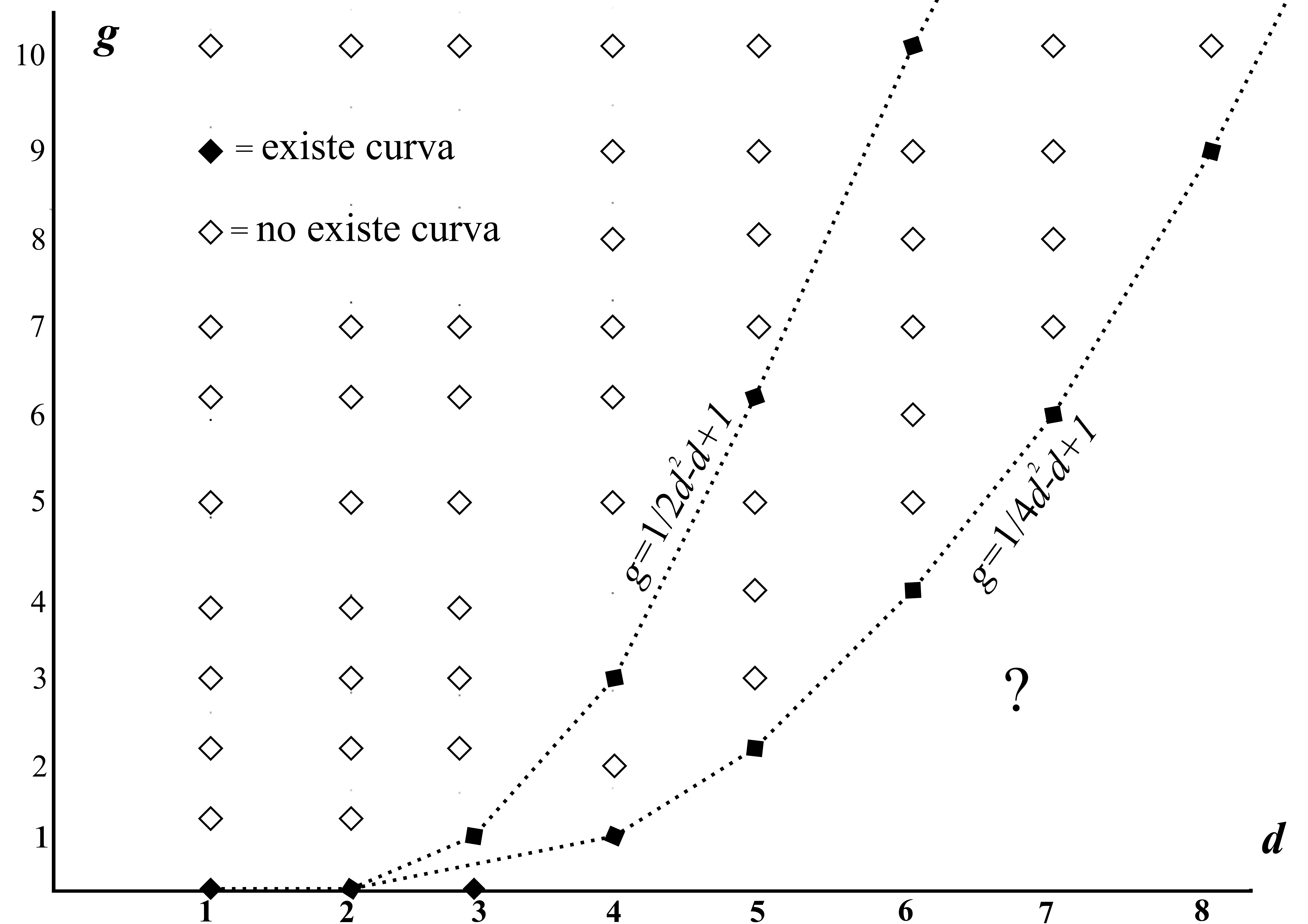}}
\caption{Curvas de grado $d$ y g\'enero $g$ en $\mathbb{P}^3$.}\label{GRAFICA}
\end{figure}
\end{center}

\section{Respuesta a la pregunta de Halphen}
\bigskip\noindent
Los argumentos expuestos hasta el momento no resuelven completamente la pregunta de Halphen pues \'estos no proporcionan criterios para saber cuando un par $(d,g)$ s\'i ocurre. Por ejemplo, estos m\'etodos son insuficientes para concluir si existe una curva suave de grado $7$ y g\'enero $5$. Sin embargo, profundizando en estas ideas en \cite{HAL} se listan todos los posibles g\'eneros de las curvas en $\PP^3$ de grado $\le 20$. La cota de Castelnuovo, pese a no resolver la pregunta de Halphen, gener\'o ideas que a\'un motivan investigaci\'on \cite{EH,CP,IVAN}. 

\medskip\noindent
La pregunta de Halphen fue finalmente resuelta por L. Gruson y C. Peskine en \cite{GP}. En este art\'iculo, se dan condiciones necesarias y suficientes para que el par de n\'umeros $(d,g)$ ocurran con el grado y g\'enero de una curva suave $C\subset \PP^3$. Dichas condiciones son: una curva suave $C\subset \PP^3$ de grado $d$ y g\'enero $g$ no contenida en ninguna superficie cu\'adrica  existe si y s\'olo si  $$0\le g\le \tfrac{1}{6}d^2-\tfrac{1}{2}d+1.$$ Esto contesta completamente la pregunta de Halphen para curvas en el espacio proyectivo $\PP^3$. Si volvemos a la figura \ref{GRAFICA}, esta soluci\'on se ve como una par\'abola adicional a las dos que aparecen, a saber $g=\tfrac{1}{6}d^2-\tfrac{1}{2}d+1$, y la regi\'on bajo \'esta es donde aparecen pr\'acticamente todas las curvas suaves de $\PP^3$.  La idea para demostrar este resultado es identificar el rango de posibles valores de $g$ que pueden ocurrir como el g\'enero de una curva suave de grado $d$ que yace en una superficie de un grado dado. Esto fue exactamente lo que hicimos en la secci\'on \ref{UNO} para curvas en un plano, y en la secci\'on \ref{SECTION2} con las curvas de Castelnuovo las cuales est\'an contenidas en superficies cu\'adricas. 

\medskip\noindent
Como dato curioso, en uno de los libros sobre geometr\'ia algebraica m\'as conocidos Robin Hartshorne aborda la pregunta de Halphen \cite[p.~349]{HART} y lista, en una gr\'afica similar a la de la figura \ref{GRAFICA}, varios pares $(d,g)$ para los cuales dice: ``no se sabe si existe una curva suave con este grado y g\'enero''. Por supuesto, \cite{HART} se public\'o en 1977 y \cite{GP} cinco a\~nos despu\'es.

\medskip\noindent
El lector observar\'a que hay cien a\~nos entre las publicaciones donde se plantea la pregunta original de Halphen en $\PP^3$ y donde \'esta se resuelve \cite{GP,HAL}. Lo que sucede para otros espacios proyectivos es actualmente foco de investigaci\'on.
Es decir, dado un $n$, nos podemos preguntar por las condiciones necesarias y suficientes para que los n\'umeros $(d,g)$ ocurran como el grado y el g\'enero de una curva suave $C\subset \PP^n$. La soluci\'on de los casos $n=2,3$ se expone en este art\'iculo trabajando sobre los n\'umeros complejos (Hartshorne resolvi\'o el caso de caracter\'istica positiva). Rathmann resolvi\'o los casos $n=4,5$, usando los m\'etodos de Gruson y Peskine. Estos son los casos donde la pregunta de Halphen tiene una respuesta completa y satisfactoria. Hasta donde tenemos conocimiento, para $n> 6$, la pregunta de Halphen sigue abierta en general. V\'ease \cite{CIRO}, donde se ataca el caso $n>6$ para ciertos valores de $(d,g)$.

\section*{Ep\'ilogo}
\medskip\noindent
Para exponer las ideas de este art\'iculo, usamos las coordenadas del espacio proyectivo $\PP^n$; lo cual era contrario a los lineamientos originales del Premio Steiner. 
Despu\'es de que todos los estudiantes y seguidores de Steiner murieron, un ge\'ometra que no hiciera uso de coordenadas era cosa dif\'icil de hallar. Por tal raz\'on, los lineamientos del Premio Steiner cambiaron diendo paso a t\'ecnicas m\'as modernas (como las de Halphen), dejando poco a poco los m\'etodos sint\'eticos en el olvido. 

\medskip\noindent
Actualmente, para estudiar curvas en $\PP^n$ se combinan aspectos tanto extr\'insecos, (usando coordenadas), como intr\'insecos (libres de coordenadas). Un ejemplo prominente de esta combinaci\'on de puntos de vista es el estudio de la geometr\'ia de ciertas curvas de Castelnuovo llamadas \textit{curvas can\'onicas} de g\'enero $g$ contenidas en $\PP^{g-1}$. El estudio de estas curvas, concentra mucha actividad e investigaci\'on en estos d\'ias en geometr\'ia algebraica. Abordar este tema ser\'ia una continuaci\'on natural de este escrito que lamentablemente no haremos. Sin embargo, el lector puede consultar \cite{ACGH,GH,MC,HS,MUM} donde se abunda en esta direcci\'on; y en muchas muchas otras.

\section*{Agradecimientos}

\medskip\noindent
Este art\'iculo se escribi\'o con el financiamiento de las C\'atedras 
\begin{footnotesize}CONACYT,\end{footnotesize} 2014-01.
Agradezco al profesor Joe Harris por introducirme en este tema. Le doy las gracias a Rosa Elisa T Hernandez Acosta de la Facultad de Ciencias de la \begin{footnotesize}UNAM\end{footnotesize} por leer cuidadosamente los borradores de este texto y, con ello, mejorar sustancialmente el espa\~nol de este art\'iculo.

\section*{Bibliograf\'ia introductoria adicional}

\noindent
Adem\'as de las refrencias introductorias incluidas en el texto, a continuaci\'on listamos 5 textos de nivel b\'asico que consideramos importantes y que no se mencionan en el art\'iculo.

\medskip\noindent
\begin{enumerate}
\item  B.~Hassett. \emph{Introduction to Algebraic Geometry.} Cambridge University Press, 2007.

\noindent
Introducci\'on a la geometr\'ia algebraica haciendo referencia siempre a problemas fundamentales del \'area y temas actuales de investigaci\'on.

\item  K.~Hulek. \emph{Elementary Algebraic Geometry.} Student Mathematical Library, AMS. 2000.

\noindent
Este es un texto introduct\'orio a la geometr\'ia algebrica desde un punto de vista moderno. Muy recomendable.

\item  F.~Kirwan. \emph{Complex Algebraic Curves.} London Mathematical Society Student Texts \text{23}, 1992.

\noindent
Este texto es una introducci\'on a la teor\'ia de curvas algebraicas planas sobre los n\'umeros complejos. En este texto se calcula, usando topolog\'ia, la f\'ormula grado-g\'enero que nosotros dedujimos en la secci\'on \ref{UNO}.

\item  R.~Miranda. \emph{Algebraic Curves and Riemann Surfaces.} Graduate Studies in Mathematics, AMS. \text{5}, 1995.

\noindent
Texto que aborda muchos aspectos de las curvas algebraicas, en particular estudia a las curvas usando an\'alisis complejo.

\item M.~Reid. \emph{Undergraduate Algebraic Geometry.} London Mathematical Society Student Texts \text{12}, 1988.

\noindent
Excelente introducci\'on a la geometr\'ia algebraica de nivel licenciatura. 
\end{enumerate}

\end{document}